\documentclass[11pt]{article}
\usepackage{amssymb}
\usepackage{amsthm}
\usepackage{graphicx}
\usepackage{graphics}
\usepackage{amsthm}

\input pictex.tex

\theoremstyle{definition}

\input epsf

\begin{document}

\date{}
\author{Andr\'es Navas}

\title{Une remarque \`a propos de l'\'equivalence bilipschitzienne entre des ensembles de Delone\\
\vspace{0.3cm}
A remark concerning bi-Lipschitz equivalence of Delone sets}
\maketitle

\begin{small}
\noindent{\bf R\'esum\'e.} Nous d\'emontrons que tout ensemble de Delone lin\'eairement r\'epetitif est rectifiable par un 
hom\'eomorphisme bilipschitzien de l'espace qui l'envoie sur l'ensemble des points \`a coordonn\'ees enti\`eres.

\vspace{0.1cm}

\noindent{\bf Abstract.} Linearly repetitive Delone sets are shown to be rectifiable by a bi-Lipschitz homeomorphism 
of the Euclidean space sending it to the standard lattice.
\end{small}

\vspace{0.5cm}

{\large\noindent{\bf Version anglaise abr\'eg\'ee}}

\vspace{0.3cm}

Over the last years, the study of {\em Delone sets} (that is, uniformly separated, uniformly discrete subsets of $\mathbb{R}^d$) has attracted the attention of many people, 
both because of their intrinsic geometric interest and their relation with mathematical models of quasi-crystrals. One of the most basic questions concerns rectifiability, 
that is, bi-Lipschitz equivalence of such a set $\mathcal{D} \subset \mathbb{R}^d$ with the standard lattice, which means the existence of a bi-Lipschitz map 
$f : \mathcal{D} \rightarrow \mathbb{Z}^d$. Examples of non-rectifiable Delone sets (for $d \geq 2$) were independently shown to exist by Burago and Kleiner 
[4] and Mc Mullen [8] (more concrete examples appear in [6]). {\em A priori}, the existence of a bi-Lipschitz homeomorphism $F : \mathbb{R}^d \to \mathbb{R}^d$ 
extending $f$ above is a stronger condition. This problem was raised (in much more generality) in [1], where it is shown that such an $F$ exists (which good 
control on the bi-Lipschitz constant) whenever the bi-Lipschitz constant of $f$ is very close to 1. However, the general case remains open. 

In this Note, we provide a quite short and elementary argument showing the existence of such an $F$ for particular choices of $f$ that are known to exist for a 
relevant class of rectifiable Delone sets (which includes, for instance, the set of vertices of Penrose tilings). Recall that a Delone set $\mathcal{D}$ is said to 
be {\em repetitive} if for each $r > 0$ there exists $R = R(r) > 0$ such that for every pair of balls $B_r,B_R$ of radii $r$ and $R$, respectively, the intersection  
$B_R \cap \mathcal{D}$ contains a translated copy of $B_r \cap \mathcal{D}$. The set $\mathcal{D}$ is said to be {\em linearly repetitive} if $R$ can 
be taken to be linear as a function of $r$. Examples of linearly repetitive Delone sets include the set of vertices of substitution tilings 
and certain cut-and-project Delone sets; see [2]. 

\vspace{0.35cm}

\noindent{\bf Main Theorem.} If $\mathcal{D}$ is a linearly repetitive Delone set in $\mathbb{R}^d$, then there exists a bi-Lipschitz homeomorphism 
$F : \mathbb{R}^d \rightarrow \mathbb{R}^d$ such that $F (\mathcal{D}) = \mathbb{Z}^d$.

\vspace{0.35cm}

This is an easy exercise for $d = 1$. For $d \geq 2$, the proof 
follows as a concatenation of results and remarks that are already known plus an elementary lemma, as next explained. 
First, Lagarias and Pleasants prove in [7] several ergodic type estimates for linearly repetitive Delone sets that imply in particular that Burago-Kleiner's 
condition from [5] holds, as shown along the proof of Theorem 2.1 in [3] (see also [9] for the case of Penrose tilings). Thus, by a theorem first proved in [5] 
for $d = 2$, and later extended to $d > 2$ in [3], such a Delone set is rectifiable. However, a careful reading of the proof in Section 4 of [5] shows that a 
stronger statement holds: if $\mathcal{D}$ satisfies the Burago-Kleiner condition, then there exists a bi-Lipschitz homeomorphism $H$ of $\mathbb{R}^d$ such 
that $H(\mathcal{D})$ is a boundedly displaced image of $\mathbb{Z}^d$. More precisely, there is a bijection $\Phi : H(\mathcal{D}) \rightarrow \mathbb{Z}^d$ 
such that $|\Phi(H(v)) - H(v)|$ is uniformly bounded on $v \in \mathcal{D}$. Knowing this, the Main Theorem follows from the next lemma applied to 
$\mathcal{D}' := \Phi (\mathcal{D})$.

\vspace{0.35cm}

\noindent{\bf Lemma.} If $\mathcal{D}'$ is a Delone set in $\mathbb{R}^d$ that is a boundedly displaced image of $\mathbb{Z}^d$, 
then there exists a bi-Lipschitz homeomorphism $F : \mathbb{R}^d \rightarrow \mathbb{R}^d$ for which $F (\mathcal{D}') = \mathbb{Z}^d$.


\vspace{0.5cm}

{\large{\noindent{\bf Version fran\c{c}aise}}

\vspace{0.3cm}

Au cours des derni\`eres ann\'ees, l'\'etude des {\em ensembles de Delone} (c'est-\`a-dire, des ensembles uniform\'ement s\'epar\'es et uniform\'ement 
discrets de $\mathbb{R}^d$) a attir\'e beaucoup d'attention en raison de leur int\'er\^et g\'eom\'etrique intrins\`eque ainsi que leur relation avec des 
mod\`eles math\'ematiques des quasicristaux. L'une des questions basiques concerne la rectifiabilit\'e, c'est-\`a-dire, l'\'equivalence d'un tel 
ensemble $\mathcal{D} \subset \mathbb{R}^d $ avec le r\'eseau standard par une application bilipschitzienne $f: \mathcal{D} \rightarrow \mathbb{Z}^d$. 
L'existence d'ensembles de Delone non-rectifiables (pour $ d \geq 2$) a \'et\'e ind\'ependamment montr\'ee par Burago et Kleiner [4] et Mc Mullen [8] 
(des exemples plus concrets apparaissent dans [6]). 
{\em \`A priori}, l'existence d'un hom\'eomorphisme bilipschitzien $F: \mathbb{R}^d  \to \mathbb{R}^d$ qui \'etend $f$ ci-dessus est une condition 
plus forte. Ce probl\`eme a \'et\'e soulev\'e (de fa\c{c}on beaucoup plus g\'en\'erale) dans [1], o\`u il est montr\'e qu'un tel $F$ existe (avec un bon 
contr\^ole de la constante bilipschitzienne) lorsque la constante bilipschitzienne de $f$ est tr\`es proche de 1. Cependant, le cas g\'en\'eral reste ouvert. 

Dans cette Note, nous donnons un argument court et \'el\'ementaire montrant l'existence d'un tel $F$ pour des choix particuliers de $f$ dont on 
conna\^it d\'ej\`a l'existence pour une classe importante d'ensembles de Delone rectifiables (qui comprend, par exemple, l'ensemble des sommets des 
pavages de Penrose). Rappelons qu'un ensemble de Delone $\mathcal{D}$ est dit {\em r\'ep\'etitif} si pour tout $r > 0$, il existe $R = R (r) > 0 $ tel que 
pour chaque paire de boules $ B_r$ et $B_R $ de rayons $ r $ et $ R $ respectivement, l'intersection $B_R \cap \mathcal{D}$ contient une copie translat\'ee 
de $B_r \cap \mathcal{D}$. L'ensemble $\mathcal{D}$ est dit {\em lin\'eairement r\'ep\'etitif} \hspace{0.04cm} si $R$ peut \^etre pris comme une fonction lin\'eaire 
de $r$. Des exemples d'ensembles de Delone lin\'eairement r\'ep\'etitifs sont l'ensemble des sommets des pavages de substitution et certains ensembles 
obtenus par coupure et projection; voir [2].

\vspace{0.35cm}

\noindent {\bf Th\'eor\`eme Principal.} Si $\mathcal{D}$ est un ensemble Delone lin\'eairement r\'ep\'etitif dans $\mathbb{R}^d$, 
alors il existe un hom\'eomorphisme bilipschitzien $F: \mathbb{R}^d \rightarrow \mathbb{R}^d $ tel que $F (\mathcal{D}) = \mathbb{Z}^d$.

\vspace{0.35cm}

Ceci est un exercice facile pour $d = 1$. Pour $d \geq 2$, la preuve r\'esulte d'une concat\'enation de r\'esultats et de remarques 
d\'ej\`a connus plus un lemme \'el\'ementaire, comme nous l'expliquons \`a continuation. Tout d'abord, Lagarias et Pleasants montrent dans [7] 
plusieurs estim\'ees de type ergodique pour les ensembles de Delone lin\'eairement r\'ep\'etitifs qui impliquent en particulier que la condition de 
Burago et Kleiner de [5] est satisfaite, comme le montre la preuve du Th\'eor\`eme 2.1 dans [3] (voir aussi [9] pour le cas des pavages de Penrose). 
Par cons\'equent, d'apr\`es un th\'eor\`eme d'abord d\'emontr\'e pour $d = 2$ dans [5], et puis \'etendu pour $ d> 2 $ dans [3], un tel ensemble de 
Delone est rectifiable. Toutefois, une lecture attentive de la preuve de la Section 4 de [5] montre qu'une conclusion plus forte est valable : si $\mathcal{D}$ 
satisfait la condition de Burago et Kleiner, alors il existe un hom\'eomorphisme bilipschitzien $H$ de $\mathbb{R}^d$ tel que $H (\mathcal{D})$ est une 
image \`a d\'eplacement born\'e de $\mathbb{Z}^d$. Plus pr\'ecis\'ement, il existe une bijection $\Phi: H (\mathcal{D}) \rightarrow \mathbb{Z}^d$ telle 
que $ | \Phi (H (v)) - H (v) |$ est uniform\'ement born\'e sur $v  \in \mathcal{D} $. Sachant cela, le Th\'eor\`eme Principal d\'ecoule 
du lemme suivant appliqu\'e \`a $ \mathcal{D}' := \Phi (\mathcal{D})$.

\vspace{0.35cm}

\noindent{\bf Lemme.} Si $ \mathcal {D}'$ est un ensemble Delone dans $ \mathbb{R}^d $ qui est une image \`a d\'eplacement born\'e de $\mathbb{Z}^d$, 
alors il existe un hom\'eomorphisme bilipschitzien $F: \mathbb{R}^d \rightarrow \mathbb{R}^d$ pour lequel $ F (\mathcal{D}') = \mathbb{Z}^d$.

\vspace {0,15 cm}

\noindent{\bf Preuve.} L'argument pour le cas $d \geq 3$ est tr\`es simple. En effet, pour chaque $\rho > 0$, le nombre de points de $\mathcal{D}'$ contenus dans toute 
boule de rayon $\rho$ est born\'e par une constante qui d\'epend de $\rho$. Par cons\'equent, si l'on bouge les points de $\mathcal{D}'$ \`a une distance uniform\'ement 
born\'ee (et tr\`es petite) par un hom\'eomorphisme bilipschitzien, on peut supposer qu'il existe $\varepsilon > 0$ de telle sorte que les segments $\ell_{v'}$ qui joignent 
les points $v'$ et $\Phi(v')$, avec $v' \in \mathcal{D}'$, sont s\'epar\'es \`a distance $\geq \varepsilon$. On peut donc consid\'erer une famille d'hom\'eomorphismes 
bilipschitziens $G_{v'}$ chacun support\'e sur le $\varepsilon / 2$-voisinage $N_{\varepsilon/2}(\ell_{v'})$ de $\ell_{v'}$ et tel que $G_{v'} (v') = \Phi (v')$. Puisque 
les longueurs des $\ell_{v'}$ sont uniform\'ement born\'ees, la constante bilipschitzienne de $G_{v'}$ peut \^etre aussi prise uniform\'ement born\'ee. Par cons\'equent, 
si l'on d\'efinit $F : \mathbb{R}^d \to \mathbb{R}^d$ en faisant $F (u) := G_{v'} (u)$ pour $u \in N_{\varepsilon / 2}(\ell_{v'} )$ et $F(u) := u$ pour $u$ 
en dehors de la r\'eunion des $N_{\varepsilon / 2} (\ell_{v'})$, nous obtenons l'hom\'eomorphisme d\'esir\'e.

\begin{figure}
  \includegraphics[width=\linewidth]{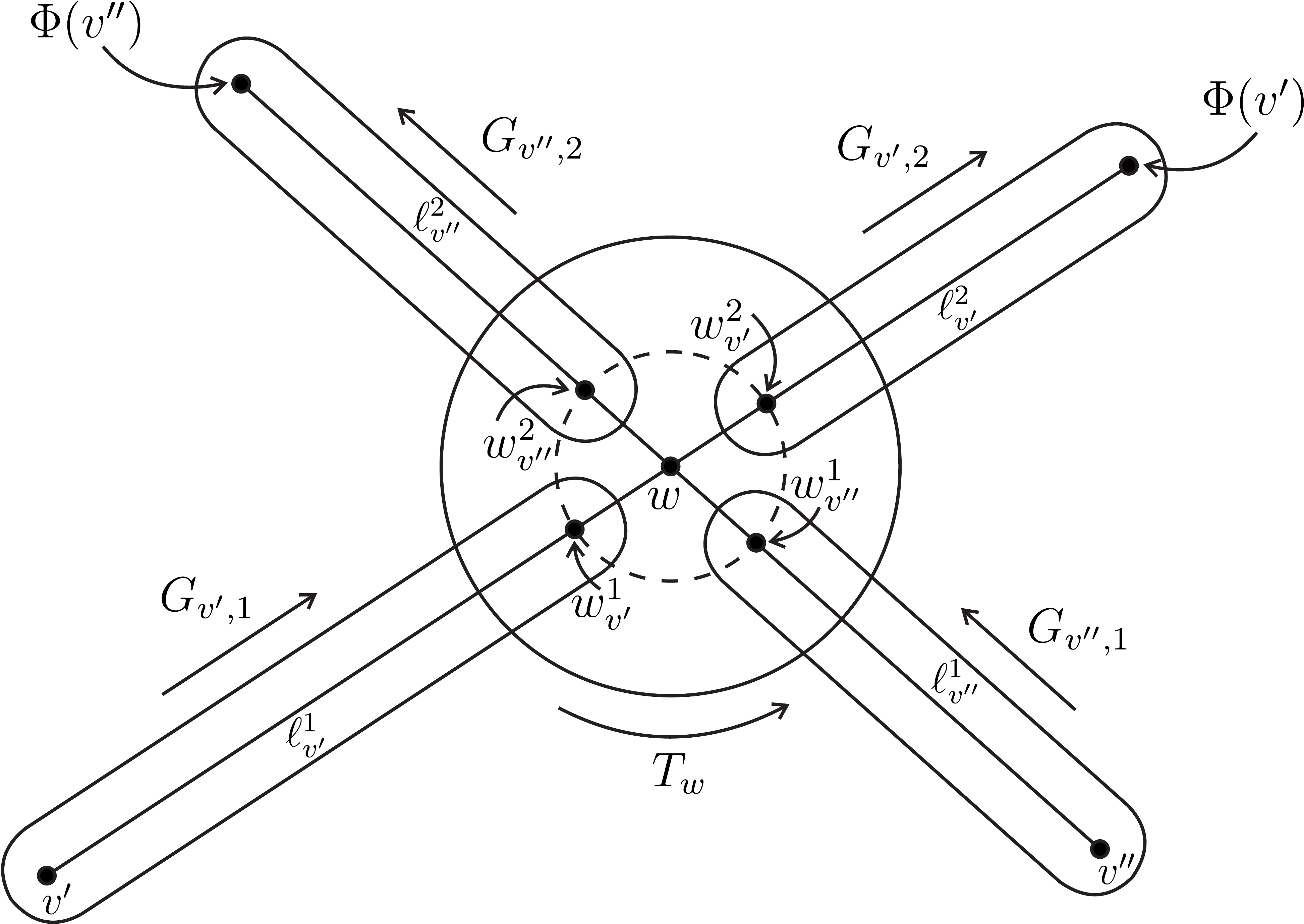}
\end{figure}

Pour le cas $d = 2$, nous ne pouvons certainement pas \'eviter les intersections des segments $\ell_{v'}$. Cependant, en faisant bouger les points de 
$\mathcal{D}'$ \`a une distance uniform\'ement born\'ee (et tr\`es petite) par un hom\'eomorphisme bilipschitzien, on peut supposer que les points 
d'intersection sont \`a distance $\geq \varepsilon$ pour un certain $\varepsilon > 0$. De plus, on peut supposer qu'ils sont \`a distance $\geq \varepsilon$ 
de tous les points $v'$ et $\Phi (v')$, et que les angles d'intersection sont aussi $\geq \varepsilon$. En supposant tout ceci, pour chaque point d'intersection 
$w$, notons $B(w)$ le disque de centre $w$ et rayon $\varepsilon / 3$. Cela donne une famille de disques disjoints, chacun desquels est le support 
d'un hom\'eomorphisme bilipschitzien $T_w$ (qui jouit d'un contr\^ole uniforme de la constante bilipschitziene) qui correspond \`a un demi-twist de Dehn,  
de telle sorte que sa restriction au cercle de centre $w$ et rayon $\varepsilon/6$ co\"{\i}ncide avec la rotation d'angle $\pi$. Notons $w^1_{v'}$ et 
$w^2_{v'}$ les points d'intersection de ce cercle et le segment $\ell_{v'}$, avec $w^1_{v'}$ entre $v'$ et $w^2_{v'}$. Remarquons que 
$T_w (w^1_{v'}) = w^2_{v'}$. Soit $\ell^i_{v'}$ le segment de $\ell_{v'}$ entre deux points cons\'ecutifs de la forme $w_{v'}^2$ et $w_{v'}^1$ 
(le pr\'emier commen\c{c}ant sur $v'$ et le dernier finissant sur $\Phi (v')$). Soit $\varepsilon' > 0$ une constante (qui ne d\'epend que de $\varepsilon$) 
telle que les $\varepsilon'$-voisinages $N_{\varepsilon'} (\ell_{v'}^{i})$ des segments $\ell_{v'}^{i}$ sont disjoints (pour tout $i$ et tout $v'$). 
De nouveau, chaque $N_{\varepsilon'} (\ell_{v'}^{i})$ est le support d'un hom\'eomorphisme bilipschitzien $G_{v',i}$ qui envoie le point initial 
de $\ell_{v'}^{i}$ sur le point final. De plus, la constante bilipschitzienne de $G_{v',i}$ peut \^etre prise uniform\'ement born\'ee. L'hom\'eomorphisme d\'esir\'e 
$F$ s'obtient alors en concatenant, dans l'ordre correct, toutes ses applications $G_{v',i}$ et $T_w$ provenant des points d'intersection $w$ le 
long de $\ell_{v'}$. Puisque le nombre de points d'intersection est uniform\'ement born\'e, cette application $F$ est un hom\'eomorphisme bilipschitzien.

\vspace {0.25cm}

\noindent{\bf Remarque}. Bien s\^ur, des calculs explicites suivant les arguments ci-dessus (en utilisant les estim\'ees d\'ej\`a connues pour l'application 
$H$) fourniraient une borne sup\'erieure pour la constante bilipschitzienne de l'hom\'eomorphisme $F$ du Th\'eor\`eme Principal en termes de la dimension 
$d \geq 2$ et de la g\'eom\'etrie de $ \mathcal{D} $ (plus pr\'ecis\'ement, des constantes de s\'eparation des points et de densit\'e relative). N\'eanmoins, 
nous ne voyons aucune application potentielle de ceci, donc nous laissons ces calculs au soin du lecteur.

\vspace {0.25cm}

\noindent{\bf Remerciements}. Je tiens \`a remercier B. Weiss pour son encouragement \`a r\'ediger cette Note, ainsi que tous les participants 
de la rencontre {\em Mathematical Quasicrystals} \`a Oberwolfach (octobre 2015). Ce travail a \'et\'e financ\'e par les projets de recherche ACT 1103, 
ACT 1415 (CONICYT) et FONDECYT 1160541.


\begin{small}

\vspace{0.15cm}

\noindent Andr\'es Navas\\ 

\noindent Dpto. de Matem\'atica y Ciencia de la Computaci\'on, Universidad de Santiago de Chile\\ 

\noindent Alameda 3363, Santiago, Chile\\ 

\noindent email: andres.navas@usach.cl

\end{small}

\end{document}